\documentclass[preprint,10pt]{elsarticle}

\usepackage{amsfonts}
\usepackage{amssymb}
\usepackage{amsmath}
\usepackage{latexsym}
\usepackage{graphicx}
\usepackage{geometry}
\usepackage{hyperref}
\usepackage{color}
\usepackage{float}
\usepackage{subfigure}
\usepackage{comment}
\usepackage{caption}

\usepackage{multicol}
\usepackage[bottom]{footmisc}
\usepackage{url}
\usepackage{bm}
\usepackage{amsthm}

\usepackage{microtype}
\usepackage{enumerate}
\usepackage[all,tips]{xy}
\SelectTips{cm}{10}
\usepackage{aliascnt}
\usepackage{epsfig}
\usepackage{epstopdf}
\usepackage{bm}
\usepackage[linesnumbered,ruled,vlined]{algorithm2e}
\usepackage{algorithmicx}
\usepackage{algpseudocode}
\usepackage{varwidth}

\newcommand{\ignore}[1]{}

\newcommand{\sfrac}[2]{\mathchoice
  {\kern0em\raise.5ex\hbox{\the\scriptfont0 #1}\kern-.15em/
   \kern-.15em\lower.25ex\hbox{\the\scriptfont0 #2}}
  {\kern0em\raise.5ex\hbox{\the\scriptfont0 #1}\kern-.15em/
   \kern-.15em\lower.25ex\hbox{\the\scriptfont0 #2}}
  {\kern0em\raise.5ex\hbox{\the\scriptscriptfont0 #1}\kern-.2em/
   \kern-.15em\lower.25ex\hbox{\the\scriptscriptfont0 #2}}
  {#1\!/#2}}

\newtheorem{my assumption}{Assumption}

\newtheorem{theorem}{Theorem}
\newtheorem{remark}{Remark}

\begin{document}

\begin{frontmatter}

\title{A Multirate Approach for Fluid-Structure Interaction Computation with Decoupled Methods}

\author{Lian Zhang\footnote{E-mail address: lzhangay@connect.ust.hk. Department of Mathematics, the HongKong University of Science and Technology.}, Mingchao Cai\footnote{E-mail address: cmchao2005@gmail.com. Department of Mathematics, Morgan State University.}
Mo Mu\footnote{E-mail address: mamu@ust.hk. Department of Mathematics, the HongKong University of Science and Technology.},
}
\begin{abstract}

We investigate a multirate time step approach applied to decoupled methods in fluid and structure interaction(FSI) computation, where two different time steps are used for fluid and structure respectively. For illustration, the multirate technique is tested by the decoupled $\beta$ scheme. Numerical experiments show that the proposed approach is stable and retains the same order accuracy as the original single time step schemes, while with much less computational expense.
\end{abstract}


\end{frontmatter}

\section{Introduction}
Fluid structure interaction (FSI) problems are extremely important because they appear in many scientific and engineering applications \cite{badia2008splitting,bazilevs2006isogeometric, bazilevs20113d,bungartz2006fluid,kuttler2008fixed,turek2006proposal, torii2006fluid}. In the literature, for solving FSI problems, both fully implicit and decoupled approaches have been applied. The fully implicit discretization approach leads to coupled schemes \cite{wick2011solving}, in which the equations of fluid dynamics, structural mechanics, and mesh moving are solved simultaneously in a fully coupled fashion. Although the coupled schemes are unconditionally stable, they result in significant difficulties and inflexibility in the design and choice of mesh generation, PDE discretization, algebraic solvers, as well as software development. In recent years, some decoupled approaches, called loosely coupled or partitioned, or explicit coupling approaches, have been developed \cite{badia2008fluid,badia2009robin,gerardo2010analysis,nobile2008effective}. In these decoupled approaches, existing fluid and structure solvers are used, the equations of fluid dynamics, structural mechanics, and mesh moving are solved sequentially or independently. However, the stability and convergence could not be guaranteed if the decoupling technique is not well-designed. For instance, in an explicit decoupling algorithm based on the Dirichlet-Neumann splitting, one solves the fluid dynamics equations with the velocity Dirichlet boundary conditions imposed by using the extrapolated value of structure velocity at the interface, then solves the structural mechanics equations with the Neumann boundary condition provided by the updated fluid interface traction, and then updates the solution of the mesh moving  by using the newest structural displacements at the interface. However, this explicit Dirichlet-Neumann scheme is known to be unconditionally unstable due to the so-called artificial added-mass effect \cite{causin2005added,forster2007artificial}. Nevertheless, note that the fluid and solid possess quite different physical properties, such as stiffness and velocity. It is natural to treat different models in their own physical regions differently for various numerical considerations. Therefore, decoupled approaches are more favorable, not only for FSI problems, but also for other coupled multi-domain, multi-physics applications \cite{cai2009numerical, mu2007, mu2010decoupled, rybak2014multiple,shan2013decoupling}.

In this work, we are interested in the coupling of an incompressible viscous fluid flow model with a thin-walled structure model. In recent years, there are two notable works: the so-called Robin-Neumann scheme and the $\beta$ scheme. In these two schemes, the interface coupling conditions are treated and approximated as a Robin type condition (a linear combination of Dirichlet condition and Neumann condition). In time marching of the Robin-Neumann scheme, one firstly solves the fluid model with the Robin interface condition approximated by using the data from the solid region at the previous time step or by certain extrapolation strategies, and then solves the structure model with the Neumann interface condition on the interface supplemented by the latest data computed from the fluid region. While, in the so-called $\beta$ scheme, the authors of \cite{bukac2016stability} split the thin-walled structure equation into two parts with a constant $\beta$. One of the two parts gives a Robin-type interface condition which is used in the fluid step while the other part can be treated as the structure equation with Neumann condition. Different from the Robin-Neumann scheme, the structure model with Neumann interface condition is solved firstly in the $\beta$ scheme. We note that, when $\beta=1$, the difference between the $\beta$ scheme and Robin-Neumann scheme exists only in that which model is solved firstly. They actually apply the same strategy for handling the interface conditions. After serous investigation and comparison, we show that the performance of the Robin-Neumann scheme is quite similar to that of $\beta$ scheme (cf. Section 4).

In this work, our main interest is to extend the $\beta$ scheme to a multirate time-stepping algorithm. By multi-rate timestepping, we mean that different time step sizes are used in different subdomains. Such a multirate time-stepping strategy is in accordance with the physical laws because the FSI problems are multi-scale problems in time. Particularly, for the Stokes flows coupled with thin-walled structures, the variables in the structure subdomain vary much more rapidly than those variables in the fluid subdomain. In the literature, a multirate time step technique was introduced in \cite{rybak2014multiple,shan2013decoupling} for coupled fluid-porous media flow models. The whole time interval $[0, T]$ is partitioned into certain coarse time grids with the time-step size $\tau_{coarse}$. Within each coarse time step, the free fluid flow solutions are computed for multiple fine time steps with the boundary information at the interface supplemented by the porous medium region (using the previous time step data). When it reaches the end of current coarse time grid, the porous medium solutions are updated by using the data from the fluid solutions. Such a multirate method is proved to be stable and convergent with the orders of accuracy in space and time depending on the spatial discretization order and time discretization order. In this work, for the Stokes flows coupled with the thin-walled structures, we choose a finer time step size for the structure model while applies a coarse time step size for the fluid flow model. Although, in a multirate time-stepping approach, one can freely choose a fine time step for either model, our numerical tests show that the current choice leads to a better numerical performance.

The paper is organized as follows. In Section 2, we describe a FSI model for coupling a Stokes flow with a thin-walled structure. In Section 3, a multirate $\beta$ scheme is outlined for the coupled FSI model. Numerical experiments are presented in Section 4 to show the stability and convergence of our scheme. Conclusions are given in Section 5.


\section{A Stokes Flow Interacting with A Thin-Walled Structure}
In this section, we describe the model problem studied in \cite{bukac2016stability, fernandez2013explicit}. In the coupled FSI model, the fluid flow motion is governed by the Stokes equations in a $d$-dimensional ($d = 2,3 $) domain $\Omega_f$ and the structure is assumed to be a linear thin-solid defined on a $(d-1)-$ manifold $\Gamma$. The boundary $\partial \Omega_f = \Gamma \cup \Gamma_D \cup \Gamma_N$ with $\Gamma_D$ and $\Gamma_N$ representing the boundaries imposed with Dirichlet and Neumann conditions respectively. The coupled model problem reads as: finding the fluid velocity $\bm{u}_f: \Omega_f \times \mathbb{R}^+ \rightarrow \mathbb{R}^\textit{d} $, the fluid pressure $ p_f: \Omega_f \times \mathbb{R}^+ \rightarrow \mathbb{R} $, and the solid displacement $ \bm{d} : \Gamma \times \mathbb{R}^+  \rightarrow \mathbb{R}^{d-1} $ such that

\begin{eqnarray}  \label{eq:fluid}
\left\{
\begin{aligned}
\rho_f \partial_t \boldsymbol{u}_f - \mbox{div} \boldsymbol{\sigma}_f  ( \boldsymbol{u_f} ,p_f ) & =0   & \text{in \quad} & \Omega_f,   \\
\textbf{div} \boldsymbol{u}_f &=  0   & \text{in \quad}       & \Omega_f ,  \\
\boldsymbol{u}_f &= 0  & \text{on \quad}  & \Gamma_D,   \\
\boldsymbol{\sigma}_f( \boldsymbol{u_f} ,p_f ) \boldsymbol{n}  &= \boldsymbol{f}_N  & \text{on \quad}  & \Gamma_N,
\end{aligned}
\right.
\end{eqnarray}
and
\begin{eqnarray}  \label{eq:solid}
\left\{
\begin{aligned}
\boldsymbol{u}_f &=\boldsymbol{u}_s = \dot{\boldsymbol{d}} & \text{on \quad} &\Gamma, \\
\rho_s \epsilon \partial_t  \dot{\bm{d}} + \boldsymbol{L}^e \bm{d} + \boldsymbol{L}^v \dot{\bm{d}}
& =- \boldsymbol{\sigma}_f( \boldsymbol{u_f} ,p_f ) \boldsymbol{n}    & \text{on \quad} &\Gamma, \\
\boldsymbol{d}  &= \boldsymbol{0}  & \text{on \quad} &\partial \Gamma,
\end{aligned}
\right.
\end{eqnarray}
satisfying the initial conditions
\[ \boldsymbol{u}_f(0) = \boldsymbol{u}_f^0, \quad \boldsymbol{d} (0) = \boldsymbol{d}^0.\]
Here, $\rho_f$ and $\rho_s$ are the fluid density and the solid density respectively, $\epsilon$ is the solid thickness,
$\boldsymbol{\dot{d}} $ is the solid velocity, $\boldsymbol{n}$ is the exterior unit normal vector to $\partial \Omega_f $,
$$
{ \varepsilon(\bm{u}_f) }  =  \frac{1}{2}  ( \nabla \bm{u}_f + \nabla \bm{u}_f^T ), \quad  \boldsymbol{\sigma}_f( \bm{u}_f, p_f) =-p_f \boldsymbol{I} + 2 \mu  \varepsilon(\bm{u}_f)
$$
with $\mu$ being the fluid dynamic viscosity, $ \boldsymbol{f}_N $ is a given surface force on $\boldsymbol{\Gamma}_N $, $ \boldsymbol{L}^e$ and $\boldsymbol{L}^v $ stand for the elastic and viscous contributions respectively. Here and hereafter, we use $\boldsymbol{L}_s = \boldsymbol{L}^e + \boldsymbol{L}^v$ to represent the solid tensor.  In the coupled model, two interface conditions are enforced: the Dirichlet condition (\ref{eq:solid})$_1$ guarantees the continuity between the fluid velocity and the structure velocity at $\Gamma$; the Neumann condition (\ref{eq:solid})$_2$ ensures the continuity of the stresses at $\Gamma$. We comment here that the equation (\ref{eq:solid})$_2$ is not only an interface coupling condition but also the structure governing equation in the coupled model.

Let ${\bm V}_f$ and ${\bm V}_s$ be the $H^1$ local spaces, associated with the appropriate global Dirichlet conditions, for the fluid and structure regions, respectively. Let $Q_f$ be the $L^2$ pressure space for Stokes model.
$$
\begin{aligned}
{\bm V}_f=\mathbf{H}_0^1 ({\Omega}_f) &= \left\{ {\bm u}_f \in (H^1({\Omega}_f))^d~|~ {\bm u}_f=0 ~\text{on}~ {\Gamma}_{D}  \right\},\\
{\bm V}_s=\mathbf{H}_0^1 ({\Omega}_s) &= \left\{ {\bm u}_s \in (H^1({\Gamma}))^{d-1}~|~ {\bm u}_s=0 ~\text{on}~ {\partial\Gamma} \right\},\\
Q_f &= L^2({\Omega_{f}}).
\end{aligned}
$$
Defining ${\bm V} \equiv \{ ({\bm u}_f, {\bm u}_s) \in {\bm V}_f \times {\bm V}_s ~|~ {\bm u}_f|_\Gamma = {\bm u}_s|_\Gamma \}$, then the implicit or monolithic weak problem for the coupled model reads as: finding $({\bm u},p_f)\in {\bm V} \times Q_f$, and ${\bm d} \in {\bm V}_s$, such that ${\bm u}_s = \frac{\partial {\bm d}}{\partial t}$ and
\begin{equation}\label{weak_form}
\begin{cases}
(\delta_t {\bm u}, {\bm v}) + a_{\Omega}({\bm u};{\bm u},{\bm d}, {\bm v}) + b({\bm v}, p_f) = f({\bm v}), ~\forall {\bm v} \in {\bm V},
\\
b({\bm u}, q) = 0, ~\forall q \in Q_f,
\end{cases}
\end{equation} where ${\bm u} \equiv ({\bm u}_f, {\bm u}_s)$, ${\bm v} \equiv ({\bm v}_f,{\bm v}_s)$,
$\delta_t {\bm u} \equiv (\rho_f  \frac{\partial {\bm u}_f}{\partial t}$, $\rho_s  \frac{\partial {\bm u}_s}{\partial t})$,
$a_{\Omega} ({\bm u};{\bm u},{\bm d}, {\bm v}) \equiv a_{\Omega_f}({\bm u}_f;{\bm u}_f,{\bm v}_f ) + a_{\Omega_s}({\bm u}_s;{\bm u}_s, {\bm d}, {\bm v}_s)$ represent the stress tensor parts, $b({\bm u},q) = - \int_{\Omega_f} q {\rm div} {\bm u}_f$.
Note that the Neumann interface condition is automatically guaranteed in the weak form, while the Dirichlet interface condition is enforced in the definition of $\bm{V}$, which reflects the coupling.

\section{Numerical Algorithms}

\subsection{The $\beta$ Scheme}
\begin{algorithm}[H]\label{beta_scheme}
	\caption{The $\beta$ scheme.}
	\begin{algorithmic}
		\State For $m=0,1,2,3...N-1$,
		\State 1. Structure step: find $\tilde{\bm u}_s^{m+1}$ such that
		\begin{equation}\label{beta_structure}
		\begin{cases}
		\rho_{s} \varepsilon \frac{\tilde{\bm u}_s^{m+1}-{\bm u}_s^m}{\Delta t }
		+\bm{L}_s ({\bm d}^{m+1})=-{\beta} {\bm \sigma}_f({\bm u}_f^{m},p_f^{m}){\bm n},&\mbox{on}~{\Gamma}, \\
		d_{t}{\bm d}^{m+1}=\tilde{\bm u}_s^{m+1},&\mbox{on}~{\Gamma}.
		\end{cases}
		\end{equation}
		\State 2. Fluid step: find ${\bm u}_f^{m+1}$, $p_{f}^{m+1}$ and ${\bm u}_s^{m+1}$ such that
		\begin{equation}\label{beta_fluid}
		\begin{cases}
		\frac{\rho_f}{\Delta t}({\bm u}_f^{m+1}-{\bm u}_f^{m})-\mbox{div}{\bm \sigma}_f({\bm u}_f^{m+1},p_{f}^{m+1})=\bm{0},&\mbox{in}~{\Omega}_f, \\
		\mbox{div}{\bm u}_f^{m+1}=0,&\mbox{in}~{\Omega}_f, \\
		{\rho}_s \varepsilon \frac{{\bm u}_s^{m+1}-\tilde{\bm u}_s^{m+1}}{\Delta t}=-{\bm \sigma}_f({\bm u}_f^{m+1},p^{m+1}_f){\bm n}+\beta{\bm \sigma}_f({{\bm u}_f^{m}, p^{m}_f}){\bm n},&\mbox{on}~{\Gamma}, \\
		{\bm u}_f^{m+1}={\bm u}_s^{m+1},&\mbox{on}~{\Gamma}. \\
		\end{cases}
		\end{equation}
	\end{algorithmic}
\end{algorithm}

In Algorithm \ref{beta_scheme}, we describe the $\beta$ scheme proposed in \cite{bukac2016stability}. The key of the $\beta$ scheme is that the structure equation is split as
\begin{eqnarray}	
\rho_s \epsilon \frac{\underbrace{{\bm u}_s^{m+1}-\tilde{\bm u}_s^{m+1}}+\tilde{\bm u}_s^{m+1}-{\bm u}_s^m}{\Delta t} +\bm{L}_s({\bm d}^{m+1}, \dot{\bm d}^{m+1})= \nonumber \\
\underbrace{-{\boldsymbol{\sigma}}_f({\bm u}_f^{m+1}, p_f^{m+1}){\bm n}+{\beta} \boldsymbol{\sigma}({\bm u}_f^m,p^n_f){\bm n}}-{\beta} \boldsymbol{\sigma}_f({\bm u}_f^n,p^n_f){\bm n}.
\end{eqnarray}
Here, the "$\underbrace{}$" parts are used in the fluid step as a Robin-type interface condition whereas the other parts are computed in the structure step.

\subsection{A Multirate $\beta$ Scheme}
\begin{figure}[H]
	\centering
	\includegraphics[scale=0.6]{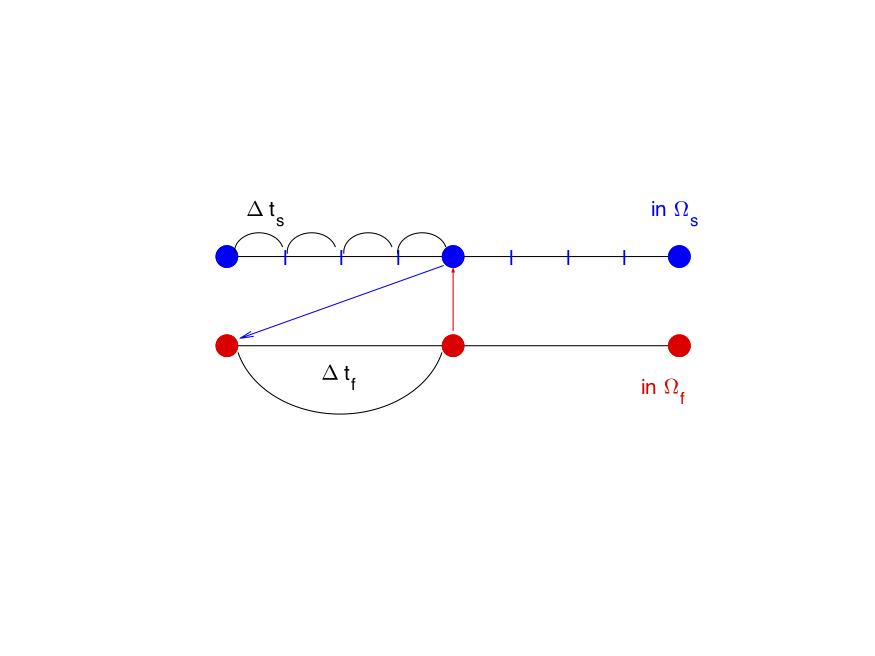}
	\caption{An illustration of a multirate time stepping technique.}
	\label{Multirate_time_step}
\end{figure}

\begin{algorithm}
\caption{A multirate $\beta$ scheme.}
\label{Multirate_beta_scheme}
\begin{algorithmic}
\State For $k=0,1,2,3...N-1$, $~~\text{set}~m_k=r\cdot k$,
\State 1. Structure steps:$~~~~\text{for}~m=m_k,m_k+1,m_k+2,...,m_{k+1}-1$,
\begin{equation}\label{Multirate_beta_structure}
	\begin{cases}
	\rho_{s} \varepsilon \frac{\tilde{\bm u}_s^{m+1}-{\bm u}_s^m}{\Delta t_s }
	+\bm{L}_s ({\bm d}^{m+1})=-{\beta} {\bm \sigma}_f({\bm u}_f^{m_k},p_f^{m_k}){\bm n},&\mbox{on}~{\Gamma}, \\
	d_{t_s}{\bm d}^{m+1}=\tilde{\bm u}_s^{m+1},&\mbox{on}~{\Gamma}.
	\end{cases}
	\end{equation}
\State 2. Fluid step:
\begin{equation}\label{Multirate_beta_fluid}
	\begin{cases}
	\frac{\rho_f}{\Delta t_f}({\bm u}_f^{m_{k+1}}-{\bm u}_f^{m_k})-\mbox{div}{\bm \sigma}_f({\bm u}_f^{m_{k+1}},p_{f}^{m_{k+1}})=\bm{0},&\mbox{in}~{\Omega}_f, \\
	\mbox{div}{\bm u}_f^{m_{k+1}}=0,&\mbox{in}~{\Omega}_f, \\
	{\rho}_s \varepsilon \frac{{\bm u}_s^{m_{k+1}}-\tilde{\bm u}_s^{m_{k+1}}}{\Delta t_f}=-{\bm \sigma}_f({\bm u}_f^{m_{k+1}},p^{m_{k+1}}_f){\bm n}+\beta{\bm \sigma}_f({{\bm u}_f^{m_{k}}, p^{m_{k}}_f}){\bm n},&\mbox{on}~{\Gamma}, \\
	{\bm u}_f^{m_{k+1}}={\bm u}_s^{m_{k+1}},&\mbox{on}~{\Gamma}. \\
	\end{cases}
	\end{equation}
\end{algorithmic}
\end{algorithm}
In the $\beta$ scheme, the coupled FSI system is split into fluid and structure steps in a sequential manner. It allows us to solve the fluid model and the structure model separately. However, in the algorithm, both the fluid solver and the structure solver use the same time step size. We note that the time scale in the structure part maybe different from the time scale in the fluid part. It is not necessary to use the same time step size in both steps. Thus, we apply a multirate time stepping technique to the $\beta$ scheme. Intuitively,
there are two possible choices of the time-stepping technique: one is to use a bigger time step size for the fluid solver, the other is using a bigger time step size for the structure  solver. Based on our numerical observations, the algorithm which uses a bigger time step size for the fluid solver whereas applies a smaller time step size for the structure solver (cf. Figure \ref{Multirate_time_step}) gives a better accuracy. We therefore describe our multirate $\beta$ scheme in Algorithm \ref{Multirate_beta_scheme} and the corresponding fully discrete weak form is given in Algorithm \ref{Multirate_beta_scheme_weakform}.

\begin{algorithm}
\caption{The fully discrete weak form for the multirate $\beta$ scheme.}
\label{Multirate_beta_scheme_weakform}
\begin{algorithmic}
\State For $k=0,1,2,3...N-1$, $~\text{set}~m_k=r\cdot k$
\State 1. Structure step: $\text{for}~m=m_k,m_k+1,m_k+2,...,m_{k+1}-1,$ find $\tilde{\bm u}_{sh}^{m+1} \in \bm V_h^s$ with $d_{t_s} \bm d_h^{m+1}=\tilde{\bm u}_{sh}^{m+1}$ such that $\forall {\bm v_{sh}} \in {\bm V_h^s}$, there holds
\begin{equation}\label{Multirate_beta_structure_weakform}
	\rho_{s} \varepsilon \left(\frac{\tilde{\bm u}_{sh}^{m+1}-{\bm u}_{sh}^m}{\Delta t_s },{\bm v_{sh}}\right)_{\Gamma}
	+a_s ({\bm d}_h^{m+1},{\bm v_{sh}})=-\beta \left({\bm \sigma}_f({\bm u}_{fh}^{m_k},p_{fh}^{m_k}){\bm n}, \bm v_{sh}\right)_{\Gamma}.
\end{equation}
\State 2. Fluid step: find $( \bm u_{fh}^{m_{k+1}}, \bm u_{sh}^{m_{k+1}}, p_{fh}^{m_{k+1}} ) \in ( \bm V_h^f, \bm V_h^s, Q_h^f )$ with $\bm u_{fh}^{m_{k+1}}\mid_{\Gamma} = \bm u_{sh}^{m_{k+1}}$ such that $\forall ( \bm v_{fh}, \bm v_{sh}, q_{fh} ) \in ( \bm V_h^f, \bm V_h^s, Q_h^f )$ with $\bm v_{fh}\mid_{\Gamma} = \bm v_{sh}$, there holds
\begin{equation}\label{Multirate_beta_fluid_weakform}
	\begin{split}
	\rho_f \left( \frac{{\bm u}_{fh}^{m_{k+1}}-{\bm u}_{fh}^{m_k}}{\Delta t_f}, \bm v_{fh} \right)_{\Omega} + a_f( \bm u_{fh}^{m_{k+1}}, \bm v_{fh})-b(p_{fh}^{m_{k+1}},\bm v_{fh})+b(q_{fh},\bm u_{fh}^{m_{k+1}})\\
	+\rho_s \varepsilon \left( \frac{{\bm u}_{sh}^{m_{k+1}}-\tilde{\bm u}_{sh}^{m_{k+1}}}{\Delta t_f}, \bm v_{sh} \right)_{\Gamma}=\beta \left( \bm \sigma_f({{\bm u}_{fh}^{m_{k}}, p^{m_{k}}_{fh}}){\bm n}, \bm v_{sh}\right)_{\Gamma}.
	\end{split}
\end{equation}
\end{algorithmic}
\end{algorithm}

\begin{remark}
If $m=m_k$, we have ${\bm u}_s^{m}={\bm u}_s^{m_k}$ from the fluid step. If $m>m_k$, we take ${\bm u}_s^{m}=\tilde{\bm u}_s^{m}$ from the structure step.
\end{remark}

\section{Numerical Experiments}
In this section, we present numerical experiments to demonstrate the convergence and stability performance of the multirate $\beta$ scheme. The benchmark test is for numerically solving a 2D pressure wave interacting with a thin-walled structure. The displacements of the interface are assumed to be infinitesimal and that the Reynolds number in the fluid is assumed to be small. The 2D fluid domain is a rectangle $\Omega_f=[0,L]\times[0,R]$ with $L=6$cm and $R=0.5$cm. The 1D structure domain is also the fluid-solid interface given by $\Gamma=[0,L]\times{R}$. See Figure 2 for the geometry configuration.
\begin{figure}[H]
	\centering
		\includegraphics[scale=0.2]{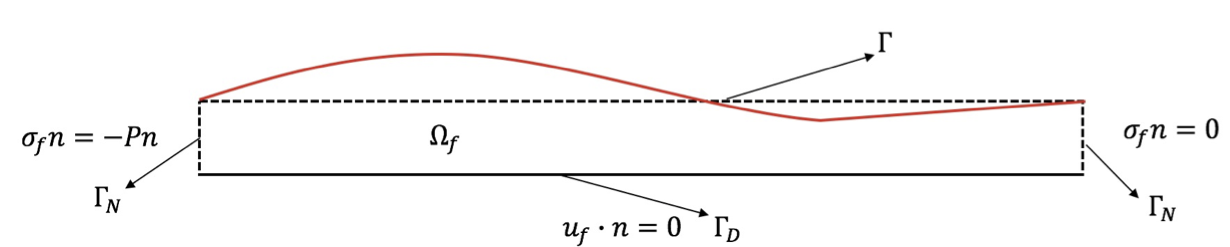}
		\caption{Geometrical configuration}
		\label{Domain}
\end{figure}

The physical parameters are: $\rho_f=1.0$, $\rho_s=1.1$, and $\mu=0.035$. The structure tensor is $\bm{L}_s({\bm d},\dot{\bm d})=c_1\partial_x^2{\bm d}+c_0{\bm d}$, $c_1=\frac{E\epsilon}{2(1+\nu)}$, $c_0=\frac{E\epsilon}{R^2(1-\nu^2)}$ with $\epsilon=0.1$, the Poisson ratio $\nu=0.5$, and the Young modulus $E=0.75\cdot10^6$. During $T^*=5 \cdot 10^{-3}$ seconds, a pressure-wave,
$$
P(t)=P_{max}(1-\cos(2t\pi/T^{*}))/2 \quad \mbox{with} \quad P_{max}=2\cdot 10^4,
$$
is prescribed on the fluid inlet boundary, a zero traction is enforced on the fluid outlet boundary, a no-slip condition is imposed on the lower boundary $y=0$. For the solid, we fix the two endpoints by imposing $\bm{d}=\bm{0}$ at $x=0$ and $x=6$.
In all the following tests, similar to \cite{fernandez2013explicit}, we generate a reference solution using the fully implicit scheme with a high space-time grid resolution $(h=3.125\times10^{-3},~\Delta t=10^{-6})$.
\begin{figure}[ht]
	\centering
	\includegraphics[scale=0.5]{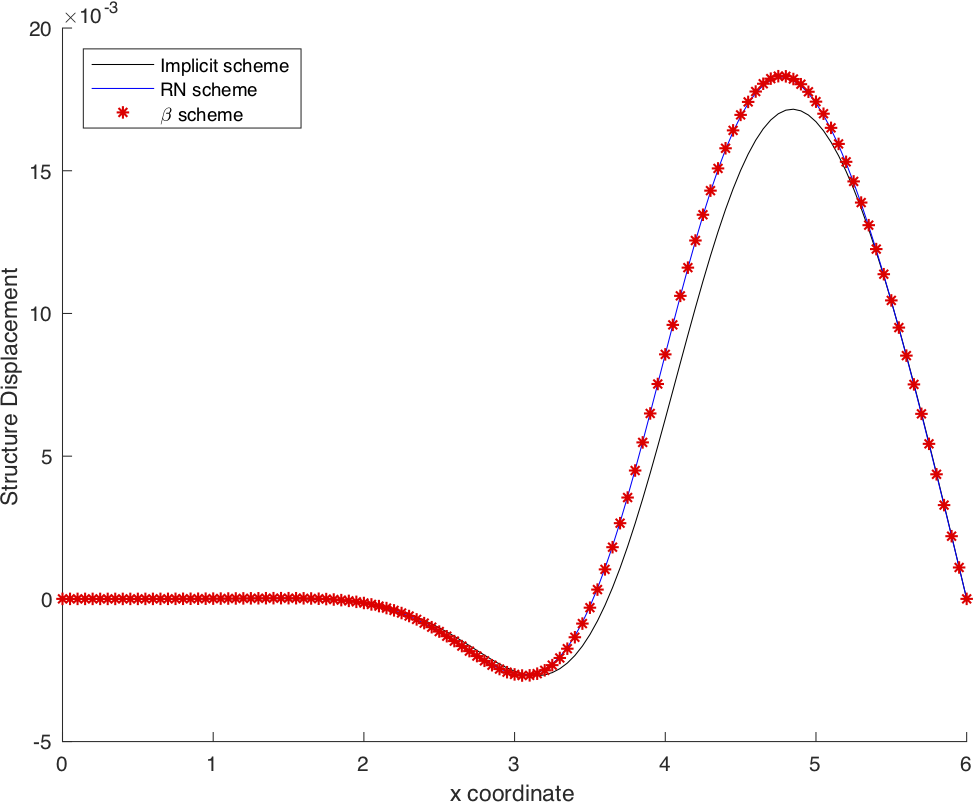}
	\caption{Comparisons of the numerical results obtained by the implicit scheme, the RN scheme and the $\beta$ scheme under the setting: $h=0.05$ and $\Delta t_s=10^{-4}$.}
	\label{Comparsion_RN_Beta}
\end{figure}

In Figure \ref{Comparsion_RN_Beta}, we compare the numerical results obtained by using the implicit scheme, the Robin-Neumann scheme and the $\beta$ scheme. The mesh size and the time step size setting for the RN and the $\beta$ scheme is: $h=0.05$ and $\Delta t=10^{-4}$. From the results, we observe that both the Robin-Neumann scheme and the $\beta$ scheme give very good approximations to the solution obtained by using the implicit scheme. Most importantly, we see clearly that the results obtained by using the Robin-Neumann scheme have little difference with those obtained by using the $\beta$ scheme. Therefore, in the following, we only report the numerical results obtained by using the $\beta$ scheme or the multirate $\beta$ scheme.

\begin{figure}[htbp]
	\centering
	\includegraphics[scale=0.4]{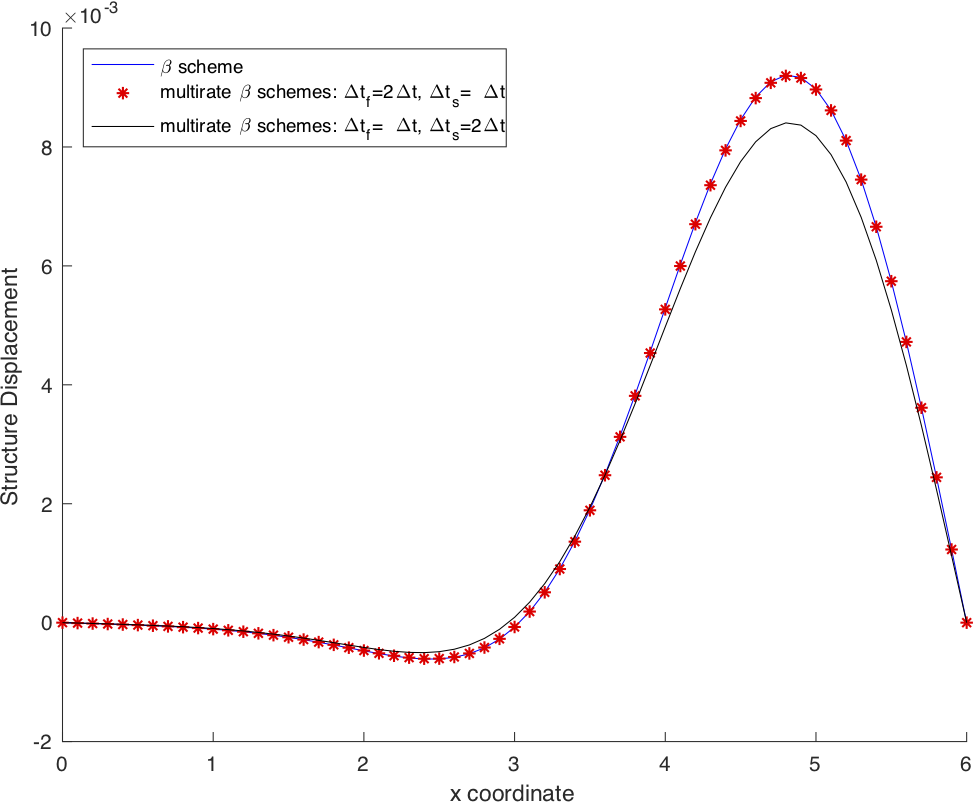}
	\caption{Comparison of the $\beta$ scheme and two different multirate $\beta$ schemes with $h=0.1, \Delta t=10^{-5}$ and $r=2$. }
	\label{Comparison_big_fluid_or_structure}
\end{figure}

In addition, we test two different multirate strategies: choosing a bigger time step size for the fluid model or a bigger time step size for the structure model. The numerical results are presented in Figure \ref{Comparison_big_fluid_or_structure}.  From the figure, we see that if the time ratio $r=2$, a bigger step size in the fluid model while applying a smaller time size for the structure model gives more accurate numerical solution than that obtained by using the other strategy. Furthermore, from our experiments and experience, the multirate $\beta$ scheme with a smaller time step size for the fluid part is unstable and the numerical  results will be messed up when $r=5$ or $10$.

\begin{figure}[h]	
	\centering
		\subfigure[]{\includegraphics[scale=0.32]{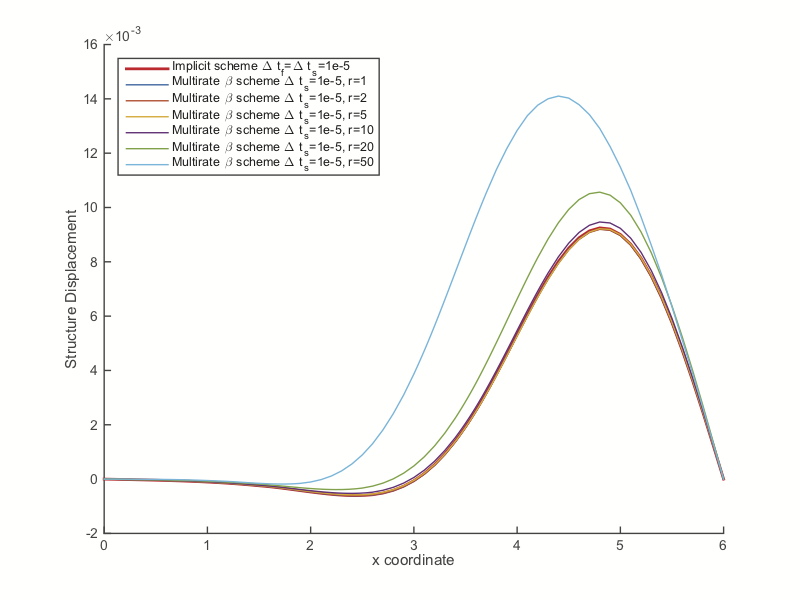}}		
		\subfigure[]{\includegraphics[scale=0.32]{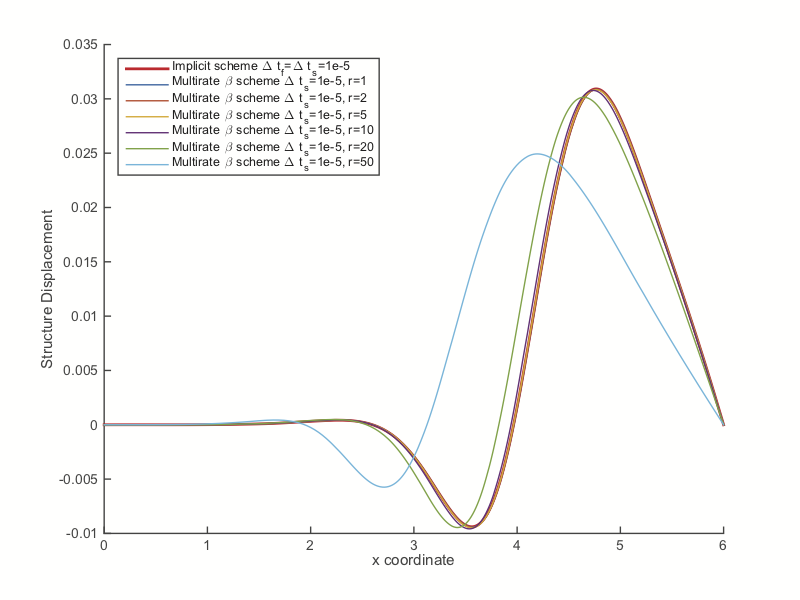}}		
		\caption{Numerical displacements under the settings: $h=0.01$ (left) $h=0.01$ (right) and $\Delta t_s=10^{-5}$.}
\label{ChangTimeStepRatio_diffh}
\end{figure}

In order to test whether a large time step ratio will cause instability, we fix $\Delta t_s$ and $h$ while vary the time ratio $r=1,5,10,20,50$. The numerical results for $t=0.015$ are reported in Figure \ref{ChangTimeStepRatio_diffh} with the structure time step size $\Delta t_s=10^{-5}$ while the mesh size $h=0.1$ (left) or $h=0.01$ (right). From the left part of the figure, we see that the solid displacement along interface obtained by the multirate $\beta$ scheme with $r=1,~2,~5,10$ are almost the same as that obtained by using the implicit scheme. Moreover, the multirate $\beta$ scheme with $r=20,~50$ are still stable although the errors become lager because of the lager time step size for the fluid model. To further investigate the stability and the convergence of the multirate $\beta$ scheme, we apply a finer mesh size $h=0.01$ (while keeping $\Delta t_s=10^{-5}$). The numerical results are presented in the right part of Figure \ref{ChangTimeStepRatio_diffh}. From the results, we have almost the same observations as those obtained under the setting $h=0.1$. Therefore, the multirate $\beta$ scheme is stable even the time size ratio is large. To have a good approximation, one only needs to keep the time step size ratio be not too large.

\begin{figure}[htbp]
	\centering
	\subfigure[$t=0.005$]{\includegraphics[scale=0.20]{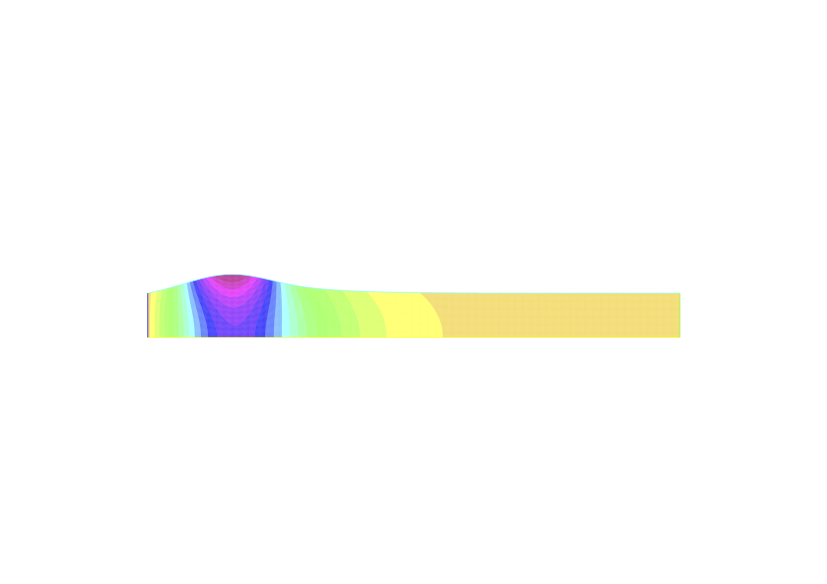}}		
	\subfigure[$t=0.010$]{\includegraphics[scale=0.20]{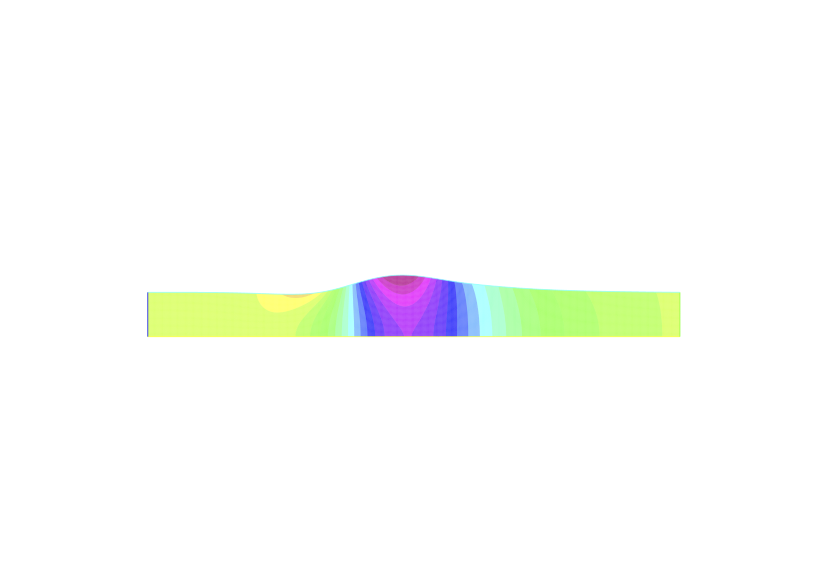}}	
	\subfigure[$t=0.015$]{\includegraphics[scale=0.20]{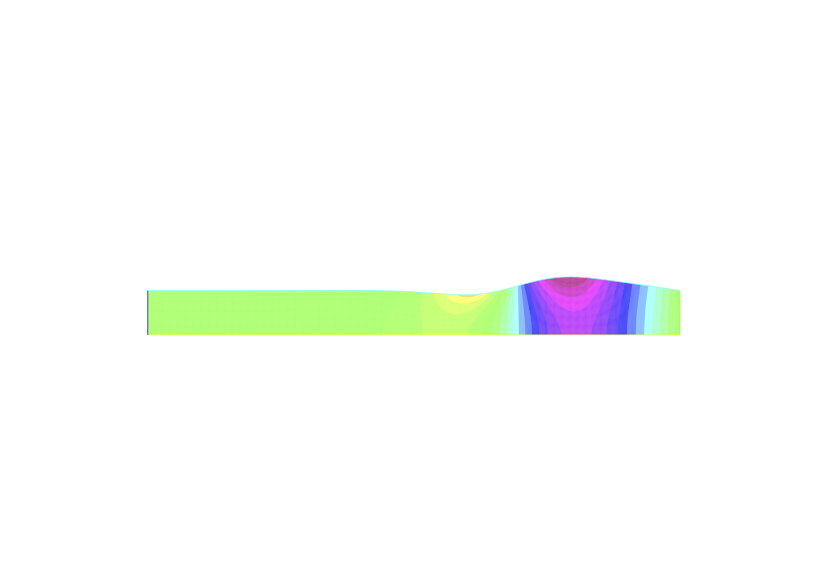}}	
	\\	
	\subfigure[$t=0.005$]{\includegraphics[scale=0.20]{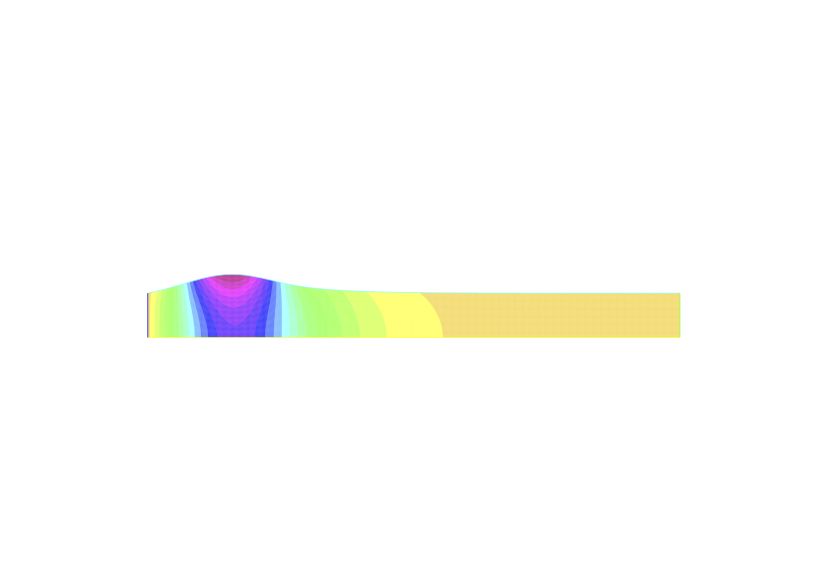}}		
	\subfigure[$t=0.010$]{\includegraphics[scale=0.20]{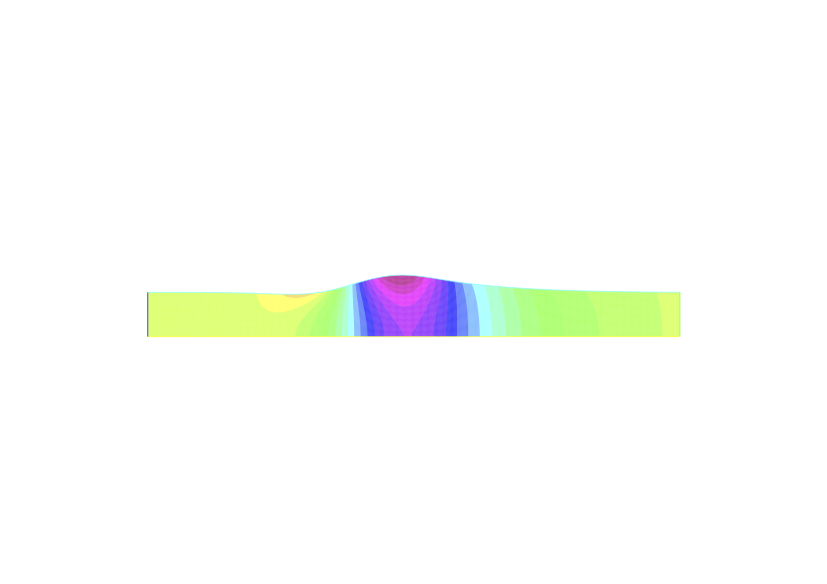}}	
	\subfigure[$t=0.015$]{\includegraphics[scale=0.20]{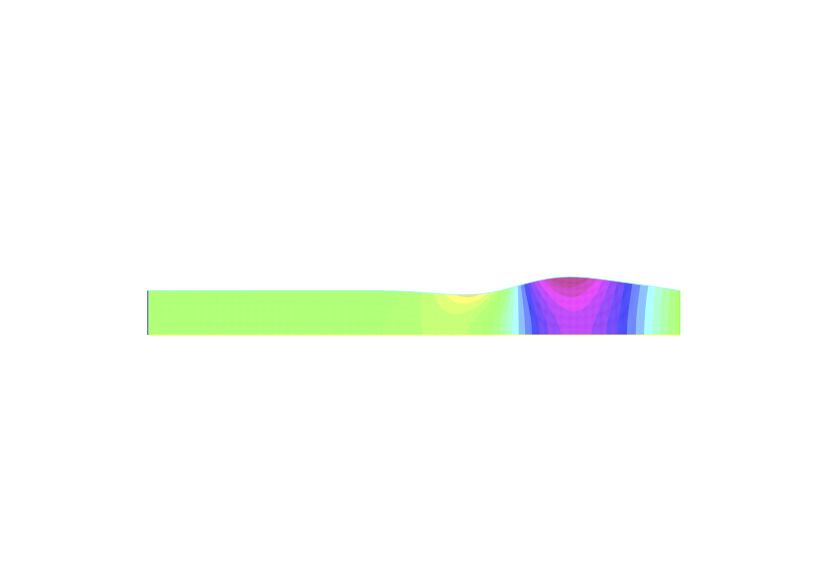}}	
	\\
	\subfigure[$t=0.005$]{\includegraphics[scale=0.20]{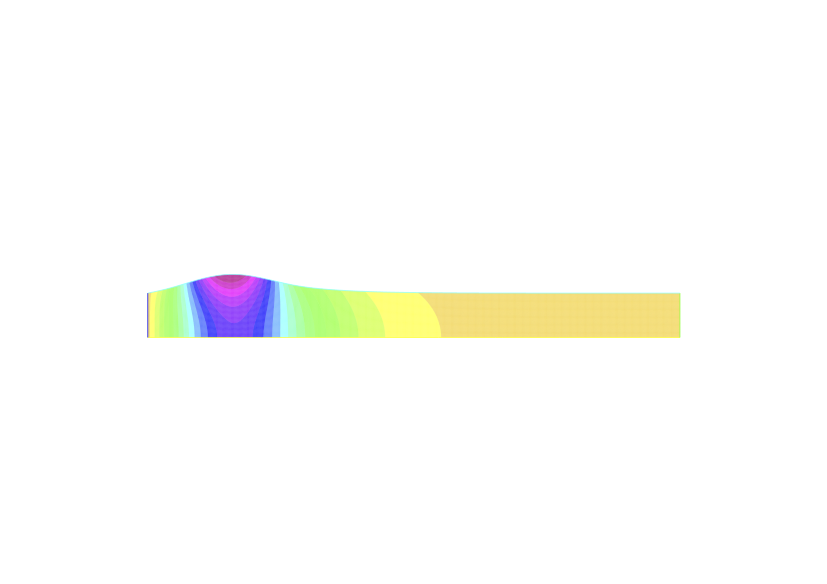}}		
	\subfigure[$t=0.010$]{\includegraphics[scale=0.20]{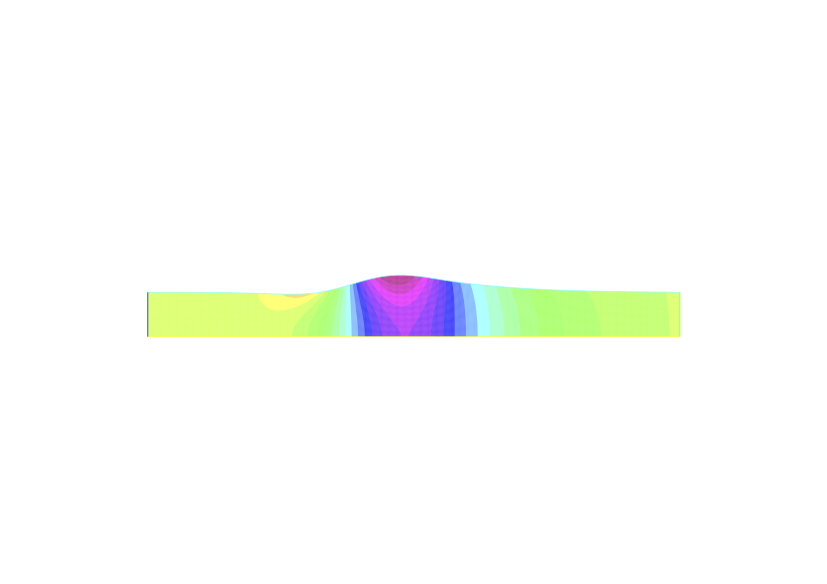}}	
	\subfigure[$t=0.015$]{\includegraphics[scale=0.20]{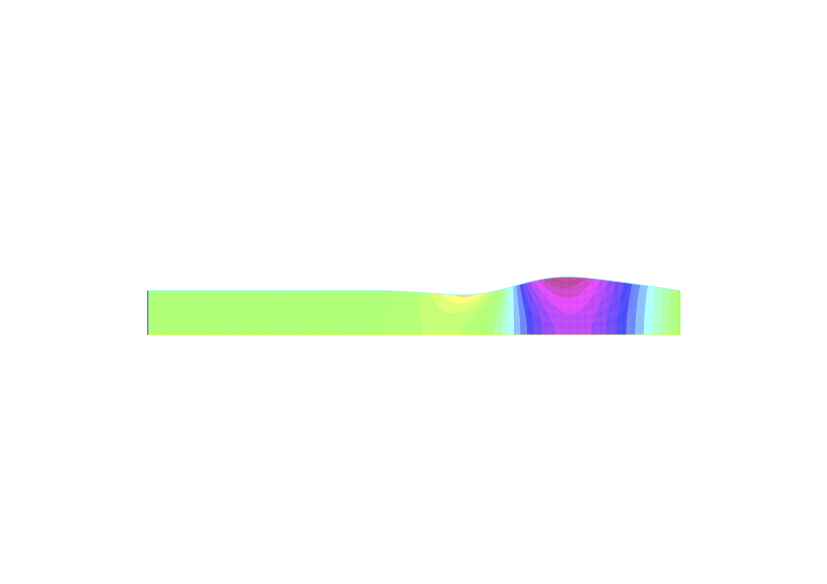}}	
	\caption{Fluid pressure distribution at $t=0.005,~ 0.010,~0.015$ obtained by the implicit scheme (top),
	the multirate $\beta$ scheme with $r=1$ (middle) and $r=10$ (bottom) with $h=0.01$ and $\Delta t_s=0.00001$.}
	\label{PressureDistribution}
\end{figure}

In Figure \ref{PressureDistribution}, for $t=0.005,~0.01,~0.015$, we compare the fluid pressure distribution obtained by using the different algorithms. From the top to the bottom, numerical results are based on the implicit scheme, the $\beta$ scheme, and the multirate $\beta$ scheme. (We comment here that the multirate $\beta$ scheme is nothing else but the $\beta$ scheme when $r=1$.) By comparing the results obtained by using different algorithms, we see that the numerical results obtained by using the multirate $\beta$ scheme are very good approximations to those obtained by using the implicit scheme.

\begin{figure}[htbp]
	\centering
		\subfigure[Relative error of $u_f$]{\includegraphics[scale=0.35]{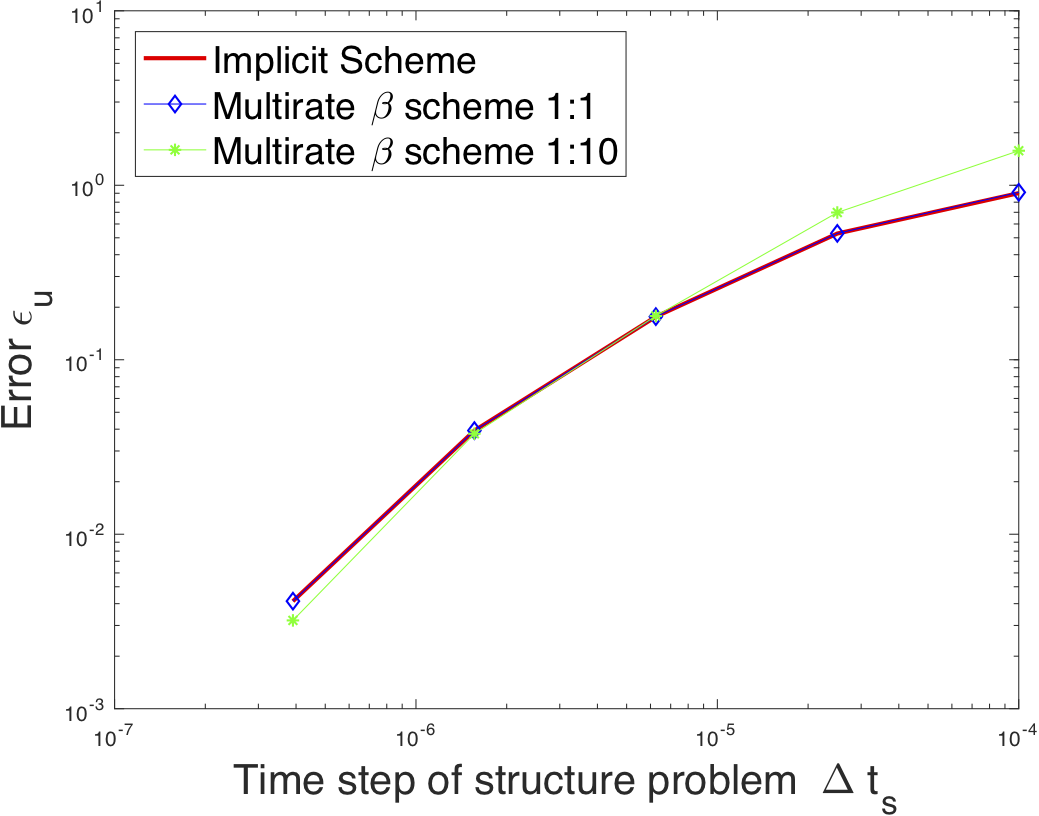}}		
		\subfigure[Relative error of $p_f$]{\includegraphics[scale=0.35]{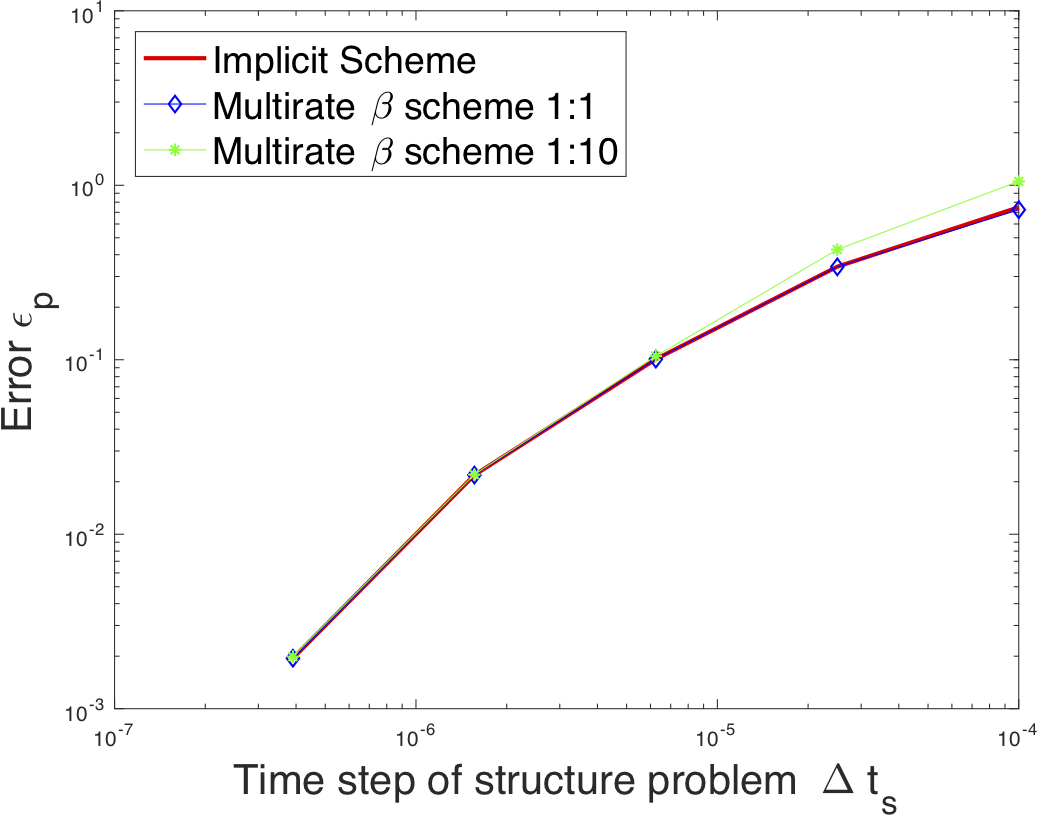}}
		\\		
		\subfigure[Relative error of $d$]{\includegraphics[scale=0.35]{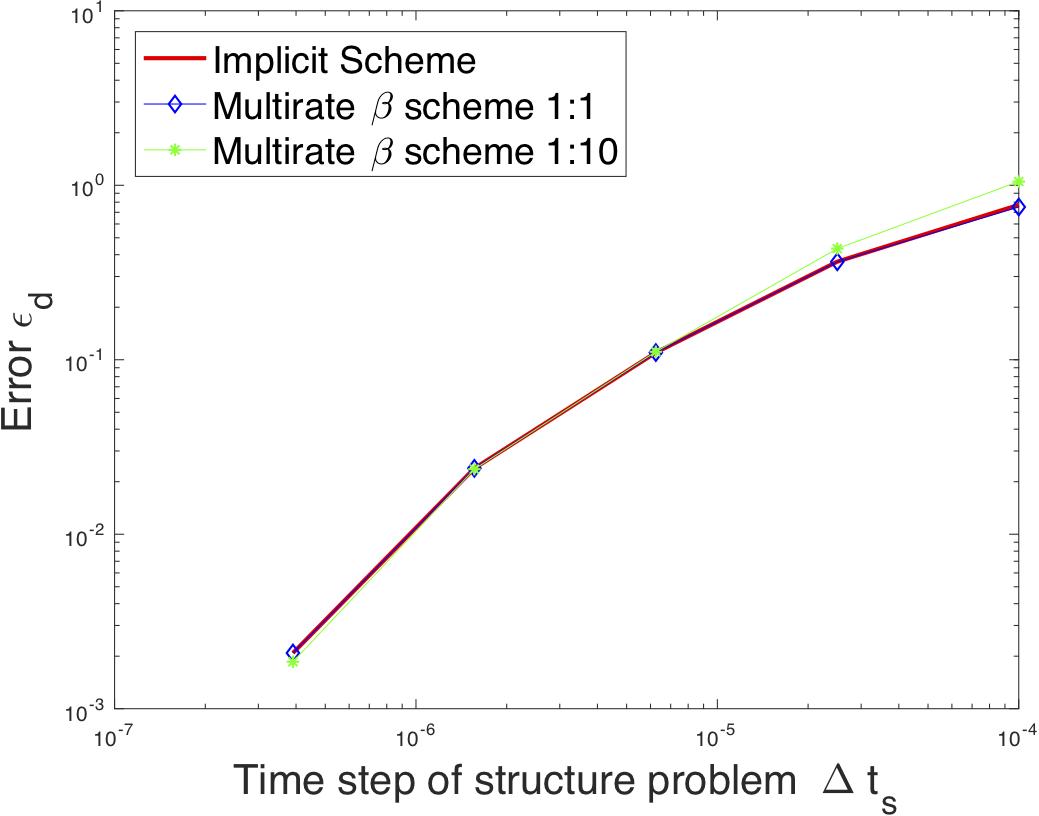}}
	\caption{Relative error of primary variables with the spacing h and time t in (\ref{refinement}).}
	\label{ChangSpacingTime}
\end{figure}

In order to examine the orders of convergence which are second order in $h$ and first order in $t$, we decrease the mesh size by a factor of two and the time step size by a factor of four at each level refinement. We start from $h=0.1$, $\Delta t_s=0.0001$ and refine four times. That is,
\begin{equation}\label{refinement}
	\{h,t_s \} =\{ 0.1 \cdot (0.5)^i , 0.0001 \cdot (0.25)^i \},~~i=0,1,2,3,4.
\end{equation}
In Figure \ref{ChangSpacingTime}, we present the relative errors of the primary variables ($\bm{u}_f$, $p_f$ and $\bm{d}$) at
$t=0.015$. We compare the reference solution (the solution obtained by the implicit scheme) with the solutions obtained by using the multirate $\beta$ scheme (with the step ratio being equal to $r=1$ or $10$). From the figures, we see that the numerical error is decreased by a factor of four as the mesh size and the time step size are refined once.

Finally, in order to highlight the advantage of the multirate $\beta$ scheme, we compare the CPU times for the different numerical algorithms under several settings of the time step sizes and mesh sizes. We fix, $\Delta t_s=10^{-5}$ and vary $h=\frac{1}{10},~\frac{1}{20},~\frac{1}{40},~\frac{1}{80},~\frac{1}{160}$. The CPU times of using different numerical algorithms are summarized in Table \ref{FixTimeChangSpaceCPU}. From the table, we see that the multirate $\beta$ scheme takes much less cost than that of the implicit scheme, in particular, when $r$ is large. Therefore, we conclude that the multirate schemes improve the efficiency.

%
%

\begin{table}[htbp]
	\caption{CPU times (in seconds) for the implicit scheme and the multirate $\beta$ scheme (with $r=1$ or $10$) under different settings of mesh sizes ($\Delta t_s=10^{-5}$ is fixed).}
	\label{FixTimeChangSpaceCPU}
	\begin{center}			
		\begin{tabular}{|c|c|c|c|}
			\hline
			&	implicit Scheme & multirate $\beta$ scheme r=1 & multirate $\beta$ scheme r=10\\
			\hline
			$h=\frac{1}{10}$    &	14.90	  & 4.02	     & 0.74		   \\
			\hline
			$h=\frac{1}{20}$	  & 48.64	  & 16.00        & 2.82		   \\
			\hline
			$h=\frac{1}{40}$	  & 179.83	  & 66.67        & 11.6		   \\
			\hline
			$h=\frac{1}{80}$	  & 797.76	  & 297.96       & 49.23		   \\
			\hline
			$h=\frac{1}{160}$	  & 3165.26	  & 1270.30      & 206.32	   \\
			\hline
		\end{tabular}	
	\end{center}
\end{table}

\section{Conclusions}
Fluid structure interaction problems appear in many engineering and science applications. Such problems are multi-domain, multi-physics problems with multiscales. In this paper, we develop a multirate $\beta$ scheme for solving the coupled model of Stokes flow interacting with a thin-walled structure. We note that the incompressible fluid model and the thin-walled structure model possess different time scales. It is natural to apply a multirate time strategy to solve such a model.
First of all, our algorithm is a decoupled algorithm which have many advantages as stated in the introduction part. Moreover, extensive numerical experiments are presented to show that the proposed scheme is efficient and accurate. Compared with the coupled implicit scheme, our algorithm uses much less computational cost to achieve the same order of accuracy.

\section{Acknowledgement}
The first and the third authors' research is supported in part by Hong Kong RGC Competitive Earmarked Research Grant HKUST16301218 and NSFC (91530319,11772281). The second author's work is supported in part by NSF Grant $\#$DMS-1831950.


\begin{thebibliography}{10}
\expandafter\ifx\csname url\endcsname\relax
  \def\url#1{\texttt{#1}}\fi
\expandafter\ifx\csname urlprefix\endcsname\relax\def\urlprefix{URL }\fi

\bibitem{badia2008fluid}
S.~Badia, F.~Nobile, C.~Vergara, Fluid--structure partitioned procedures based
  on robin transmission conditions, J. Comput. Phys. 227~(14)
  (2008) 7027--7051.

\bibitem{badia2009robin}
S.~Badia, F.~Nobile, C.~Vergara, Robin--robin preconditioned krylov methods for
  fluid--structure interaction problems, Comput. Methods Appl. Mech.
  Eng. 198~(33-36) (2009) 2768--2784.

\bibitem{badia2008splitting}
S.~Badia, A.~Quaini, A.~Quarteroni, Splitting methods based on algebraic
  factorization for fluid-structure interaction, SIAM J. Sci. Comput.30~(4) (2008) 1778--1805.

\bibitem{bazilevs2006isogeometric}
Y.~Bazilevs, V.~M. Calo, Y.~Zhang, T.~J. Hughes, Isogeometric fluid--structure
  interaction analysis with applications to arterial blood flow, Comput.
  Mech. 38~(4-5) (2006) 310--322.

\bibitem{bazilevs20113d}
Y.~Bazilevs, M.-C. Hsu, J.~Kiendl, R.~W{\"u}chner, K.-U. Bletzinger, 3d
  simulation of wind turbine rotors at full scale. part ii: Fluid--structure
  interaction modeling with composite blades, Int. J.  Numer. Methods Fluids
  65~(1-3) (2011) 236--253.

\bibitem{bukac2016stability}
M.~Bukac, B.~Muha, Stability and convergence analysis of the extensions of the
  kinematically coupled scheme for the fluid-structure
  interaction, SIAM J. Numer. Anal. 54~(5) (2016) 3032--3061.

\bibitem{bungartz2006fluid}
H.-J. Bungartz, M.~Sch{\"a}fer, Fluid-structure interaction: modelling,
  simulation, optimisation, vol.~53, Springer Science \& Business Media, 2006.


\bibitem{cai2009numerical}
M. Cai, M. Mu, and J. Xu, Numerical solution to a mixed Navier-Stokes/Darcy model by the
two-grid approach. SIAM J. Numer. Anal. 47(5) (2009) 3325-3338.

\bibitem{cai2018some}
M. Cai, P. Huang, and M. Mu, Some multilevel decoupled algorithms for a mixed
Navier-Stokes/Darcy model. Adv. Comput. Math., 44(1), (2018) pp.115-145.

\bibitem{causin2005added}
P.~Causin, J.-F. Gerbeau, F.~Nobile, Added-mass effect in the design of
  partitioned algorithms for fluid--structure problems,
  Comput. Methods Appl. Mech. Eng., 194~(42-44) (2005) 4506--4527.

\bibitem{fernandez2013explicit}
M.~A. Fern{\'a}ndez, J.~Mullaert, M.~Vidrascu, Explicit robin--neumann schemes
  for the coupling of incompressible fluids with thin-walled structures,
  Comput. Methods Appl. Mech. Eng., 267 (2013) 566--593.

\bibitem{forster2007artificial}
C.~F{\"o}rster, W.~A. Wall, E.~Ramm, Artificial added mass instabilities in
  sequential staggered coupling of nonlinear structures and incompressible
  viscous flows, Comput. Methods Appl. Mech. Eng., 196~(7)
  (2007) 1278--1293.

\bibitem{gerardo2010analysis}
L.~Gerardo-Giorda, F.~Nobile, C.~Vergara, Analysis and optimization of
  robin--robin partitioned procedures in fluid-structure interaction problems,
  SIAM J. Numer. Anal. 48~(6) (2010) 2091--2116.

\bibitem{kuttler2008fixed}
U.~K{\"u}ttler, W.~A. Wall, Fixed-point fluid--structure interaction solvers
  with dynamic relaxation, Comput. Mech. 43~(1) (2008) 61--72.


\bibitem{mu2007}
M. Mu and J. Xu, A two-grid method of a mixed Stokes/Darcy model for coupling fluid
flow with porous media flow. SIAM J. Numer. Anal., 45(5), (2007) pp.1801-1813.

\bibitem{mu2010decoupled}
M. Mu and X. Zhu, Decoupled schemes for a non-stationary mixed Stokes-Darcy model,
 Math. Comput., 79(270), (2010) 707-731.


\bibitem{nobile2008effective}
F.~Nobile, C.~Vergara, An effective fluid-structure interaction formulation for
  vascular dynamics by generalized robin conditions,
  SIAM J. Sci. Comput. 30~(2) (2008) 731--763.

\bibitem{rybak2014multiple}
I.~Rybak, J.~Magiera, A multiple-time-step technique for coupled free flow and
  porous medium systems, J. Comput. Phys. 272 (2014) 327--342.

\bibitem{shan2013decoupling}
L.~Shan, H.~Zheng, W.~J. Layton, A decoupling method with different subdomain
  time steps for the nonstationary stokes--darcy model, Numerical Methods for
  Partial Differential Equations 29~(2) (2013) 549--583.

\bibitem{torii2006fluid}
R.~Torii, M.~Oshima, T.~Kobayashi, K.~Takagi, T.~E. Tezduyar, Fluid--structure
  interaction modeling of aneurysmal conditions with high and normal blood
  pressures, Comput. Mech. 38~(4-5) (2006) 482--490.

\bibitem{turek2006proposal}
S.~Turek, J.~Hron, Proposal for numerical benchmarking of fluid-structure
  interaction between an elastic object and laminar incompressible flow, in:
  Fluid-structure interaction, Springer, 2006, pp. 371--385.

\bibitem{wick2011solving}
T.~Wick, Solving monolithic fluid-structure interaction problems in arbitrary
  {L}agrangian {E}ulerian coordinates with the deal. ii library.

\end{thebibliography}
\end{document}